\documentclass[11pt]{amsart}

\usepackage{amsmath}
\usepackage{amscd}
\usepackage{amsfonts}
\usepackage{graphicx}
\usepackage[usenames,dvipsnames]{color}
\usepackage{color}

\newtheorem{theorem}{Theorem}[section]

\newtheorem{lemma}{Lemma}[section]

\newtheorem{definition}{Definition}[section]

\newtheorem{proposition}{Proposition}[section]

\newtheorem{corollary}{Corollary}[section]

\def \epsilon {\varepsilon}
\def \phi {\varphi}
\def \proof {{\bf Proof}\hspace{.4cm}}
\def \dist {{\rm dist}}
\def \diams {{\rm diam}^\#}

\def \Aa {(A,a)}
\def \Uu {(U,u)}
\def \Upup {(U',u')}

\def \Aaaa {(A_\alpha,a_\alpha)}
\def \Uaua {(U_\alpha,u_\alpha)}

\def \ALaa {(A^{\Lambda_\alpha},a_\alpha)}

\def \Uaua {(U_\alpha,u_\alpha)}
\def \Uauap {(U'_\alpha,u'_\alpha)}
\def \Umum {(U_m,u_m)}
\def \Um {U_m}

\def \Umk {U_{m_k}}

\def \Vava {(V_\alpha,v_\alpha)}

\def \eia {\eta^i_\alpha}

\def \gia {\gamma^i_\alpha}

\def \Kia {K^i_\alpha}

\def \Kimk {K^i_{m_k}}

\def \calA {{\mathcal A}}

\def \calD {{\mathcal D}}
\def \calF {{\mathcal F}}
\def \calG {{\mathcal G}}
\def \calL {{\mathcal L}}

\def \calU {{\mathcal U}}
\def \calV {{\mathcal V}}
\def \calW {{\mathcal W}}

\def \calUp {{\mathcal U}'}
\def \calVp {{\mathcal V}'}

\def \normA {{\parallel \hspace{-.1cm}{\mathcal A} \hspace{-.1cm}\parallel}}
\def \normU {{\parallel \hspace{-.045cm}{\mathcal U} \hspace{-.09cm}\parallel}}
\def \normV {{\parallel \hspace{-.075cm}{\mathcal V} \hspace{-.075cm}\parallel}}

\def \normUE {{\parallel \hspace{-.045cm}{\mathcal U} \hspace{-.09cm}\parallel_{\mathrm E}}}

\def \Bn {{\mathfrak B}_n}
\def \Dn {{\mathfrak D}_n}
\def \Mn {{\mathfrak M}_n}
\def \inn {{\mathfrak i}_n}

\def \Car {Carath\'eodory }

\def \pt {{\rm pt}}

\def \C {\mathbb C}
\def \D {\mathbb D}
\def \DO {(\mathbb D,0)}
\def \R {\mathbb R}
\def \T {\mathbb T}

\def \cbar {\overline {\mathbb C}} 
\def \dsharpz {|{\mathrm d}^{\scriptscriptstyle \#}\hspace{-.06cm} z|}

\title{The Carath\'eodory Topology for Multiply Connected Domains II}

\author{Mark Comerford}

\address{Department of Mathematics,
University of Rhode Island,
5 Lippitt Road, Room 102F,
Kingston, RI 02881, USA.
email: {\tt mcomerford@math.uri.edu}}



\begin{document}

\renewcommand{\theenumi}{\emph{\arabic{enumi}.}}
\renewcommand{\labelenumi}{\theenumi}

\maketitle

\begin{abstract}
We continue our exposition concerning the \Car topology for multiply connected domains by introducing the notion of boundedness for a family of pointed domains of the same connectivity. The limit of a convergent sequence of $n$-connected domains which is bounded in this sense is again $n$-connected and will satisfy the same bounds. We prove a result which establishes several equivalent conditions for boundedness. 
This allows us to extend the notions of convergence and equicontinuity to families of functions defined on varying domains. \end{abstract}

\section{Bounded Families}

We start with the definition of a bounded family of non-degenerate pointed domains of the same connectivity. Essentially, a bounded family is stable with respect to limits in the \Car topology.

\begin{definition} Let $n \ge 1$ and let $\calU = \{\Uaua\}_{\alpha \in A}$ be a family of pointed domains where every domain $U_\alpha$ is an $n$-connected non-degenerate subdomain of $\cbar$. We say that  
$\calU$ is {\rm bounded} if every sequence in $\calU$ which is convergent in the \Car topology has a limit which is a non-degenerate $n$-connected pointed domain and we write $\pt \sqsubset \calU \sqsubset \cbar$. Otherwise we say that $\calU$ is {\rm unbounded}.
\end{definition}

We observe that using the Hausdorff version of \Car convergence, it is clear that any family of pointed domains (whose connectivities in general can differ) is precompact in the 
sense that any sequence will have a convergent subsequence. The content of the definition, then, is that any limit is non-trivial and non-degenerate and that connectivity is preserved. 

\scalebox{1.069}{\includegraphics{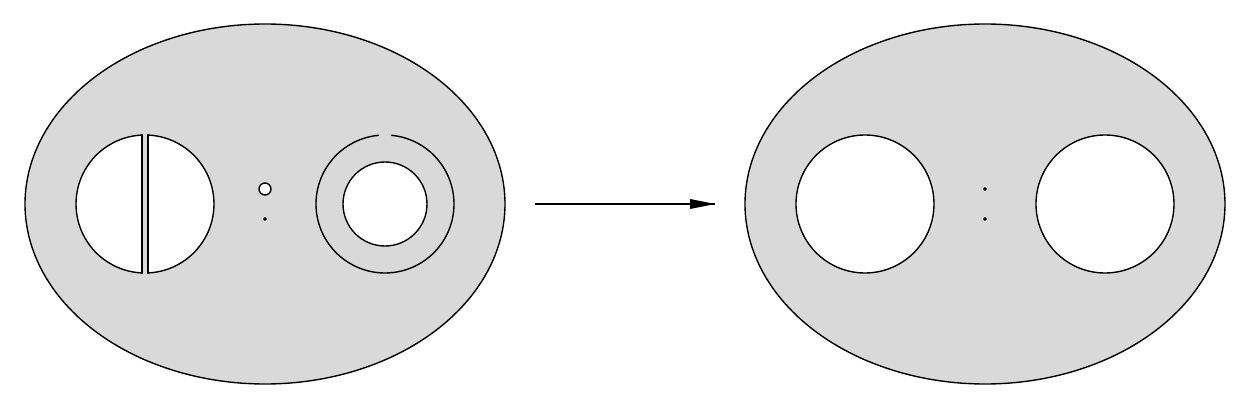}}
\unitlength1cm
\begin{picture}(0.01,0.01)
  \put(-11.5,-.3){\footnotesize$(U_m,u_m)$  }
  \put(-7.45,2.4){\footnotesize $m \to \infty$}
  \put(-3.35,-.3){\footnotesize$(U,u)$  }
  \put(-10.73,1.93){$\scriptstyle u_m$}
  \put(-2.9,1.93){$\scriptstyle u$}
\end{picture}

\vspace{.4cm}
Figure 4 above (which is the same as Figure 2 in the first part of this paper) shows an example of what can happen if a family is unbounded. 
Here we have a sequence of non-degenerate pointed domains of connectivity $6$ whose limit is a degenerate pointed domain of connectivity $4$.

It is desirable to be able to specify the boundedness of a family of pointed domains of the same connectivity in a quantitative way, the idea being that the larger the bound, the closer the members of the family are allowed to get to becoming degenerate (which includes having lower connectivity). This leads to the notion of the \emph{\Car bound} which we first define for an individual pointed domain and then for a family. Theorem 4.2 which is the main result of this part of the paper gives several equivalent conditions for boundedness, including finiteness of the \Car bound.

As we will see, it turns out that one only needs to consider the principal meridians in order to ensure bounded behaviour. However, sometimes one is interested in the other meridians as well, which leads us to also define the \emph{extended \Car bound}. We remind the reader of the convention that for a pointed domain $\Uu$ of connectivity $n$ with an extended system $\Gamma = \{ \gamma^i, 1 \le i \le E(n)\}$ of $E(n)$ meridians, we always number the meridians of $\Gamma$ 
so that the first $P(n)$ are the principal ones.

\begin{definition} Let $\Uu$ be a hyperbolic pointed domain and let $\Gamma =\{\gamma^1, \gamma^2, \ldots \ldots, \gamma^i\}$ be a collection of curves in $U$. We define the \rm{length} and \rm{distance}, $\calL(\Gamma)$ and $\calD(\Gamma)$ of $\Gamma$ by 
\vspace{.2cm}
\[ \calL(\Gamma) = \max_{\scriptscriptstyle 1 \le j \le i} |\log \ell(\gamma^j) |, \qquad 
\calD(\Gamma) = \max_{\scriptscriptstyle 1 \le j \le i} \rho(u, \gamma^j).\]

\end{definition}

\begin{definition} Let $\Uu$ be an $n$-connected non-degenerate pointed domain.  If $n \ge 2$, let $\Gamma = \{ \gamma^i, 1 \le i \le E(n)\}$ be an extended system of meridians for $U$ with lengths $l^i$ and distances $d^i$, $1 \le i \le E(n)$. We define the {\rm length} and {\rm extended length}, $\calL(U)$, $\calL_{\mathrm E} (U)$ of $U$ by 

\[ \calL(\Uu) = \max_{\scriptscriptstyle 1 \le i \le P(n)} |\log l^i |, \qquad \calL_{\mathrm E}(\Uu) = \max_{\scriptscriptstyle 1 \le i \le E(n)} |\log l^i | \]

and the {\rm distance} $\calD(\Uu)$ of $\Uu$ and {\rm extended distance} $\calD(\Uu)$, $\calD_{\mathrm E}(\Gamma)$ of $\Gamma$  by

\[ \calD(\Uu) = \max_{\scriptscriptstyle1 \le i \le P(n)} d^i, \qquad \calD_E(\Gamma) = \max_{\scriptscriptstyle1 \le i \le E(n)} d^i.\]

If $n=1$ and $U$ is simply connected, we define each of the quantities $\calL(U)$, $\calL_{\mathrm E} (U)$, $\calD(\Uu)$, $\calD_{\mathrm E}(\Gamma)$ to be zero.
\end{definition} 

Note that by Theorem 1.7 on the uniqueness of principal meridians the length and distance can be defined for $\Uu$ rather than just $\Gamma$. Also, as by Theorem 1.6 the lengths do not depend on the choice of system of meridians, we can say the same about the extended length of $\Uu$. 

The extended distances will in general depend on the choice of system, so it is of interest whether there is a system for which these distances are as small as possible. We postpone the proof of the following result until after the statement of Theorem 4.2.

\begin{lemma} For a pointed domain $\Uu$ of finite connectivity $n \ge 2$, there exists an extended system $\Gamma = \{ \gamma^i, 1 \le i \le E(n)\}$ of meridians for which the distances 
$d^i$, $1 \le i \le E(n)$ are as small as possible. 
\end{lemma}

\begin{definition} Let $\Uu$ be a pointed $n$-connected domain with $n \ge 2$. 
An extended system of meridians $\Gamma = \{ \gamma^i, 1 \le i \le E(n)\}$ as above is called {\rm maximally close}. 
\end{definition}

\begin{definition}
Let $\Uu$ be a pointed $n$-connected domain with $n \ge 1$. For $n \ge 2$ we 
define the {\rm extended distance} $\calD_{\mathrm E}(\Uu)$ of $\Uu$  by
\vspace{.25cm}
\[ \calD_{\mathrm E}(\Uu) = \max_{\scriptscriptstyle1 \le i \le E(n)} d^i .\]

where $d^i$, $1 \le i \le E(n)$ are the distances of a maximally close system $\Gamma$.

For $n=1$ where $U$ is simply connected, we define $\calD_{\mathrm E}(\Uu)$ to be zero.
\end{definition}

We observe here that using a different system other than a maximally close one to calculate $\calD_E(\Uu)$ will not give an answer that is very different. If $\gamma$ and $\gamma'$ are two meridians which are in the same homology class but different homotopy classes, then by \cite{Com2} Lemma 2.3, these two curves must intersect. If $\gamma$ and $\gamma'$ 
have lengths $l = l'$ and distances $d$ and $d'$ respectively, then it is a simple calculation to check that 
\vspace{.2cm}
\[ d' \le d + l/2.\]

Thus if $\Gamma$ and $\Gamma'$ are two systems of meridians for $\Uu$, we have
\vspace{.2cm}
\[\calD_{\mathrm E}(\Gamma') \le \calD_{\mathrm E}(\Gamma) + \frac {e^{\calL_{\mathrm E}(\Gamma)}}{2}. \]

The other remark worth making here is that even a maximally close system of meridians is not in general unique. This can be seen in Figure 5 above where both meridians are equally close to the base point but in different homotopy classes.

\vspace{.2cm}
\scalebox{0.9}{\includegraphics{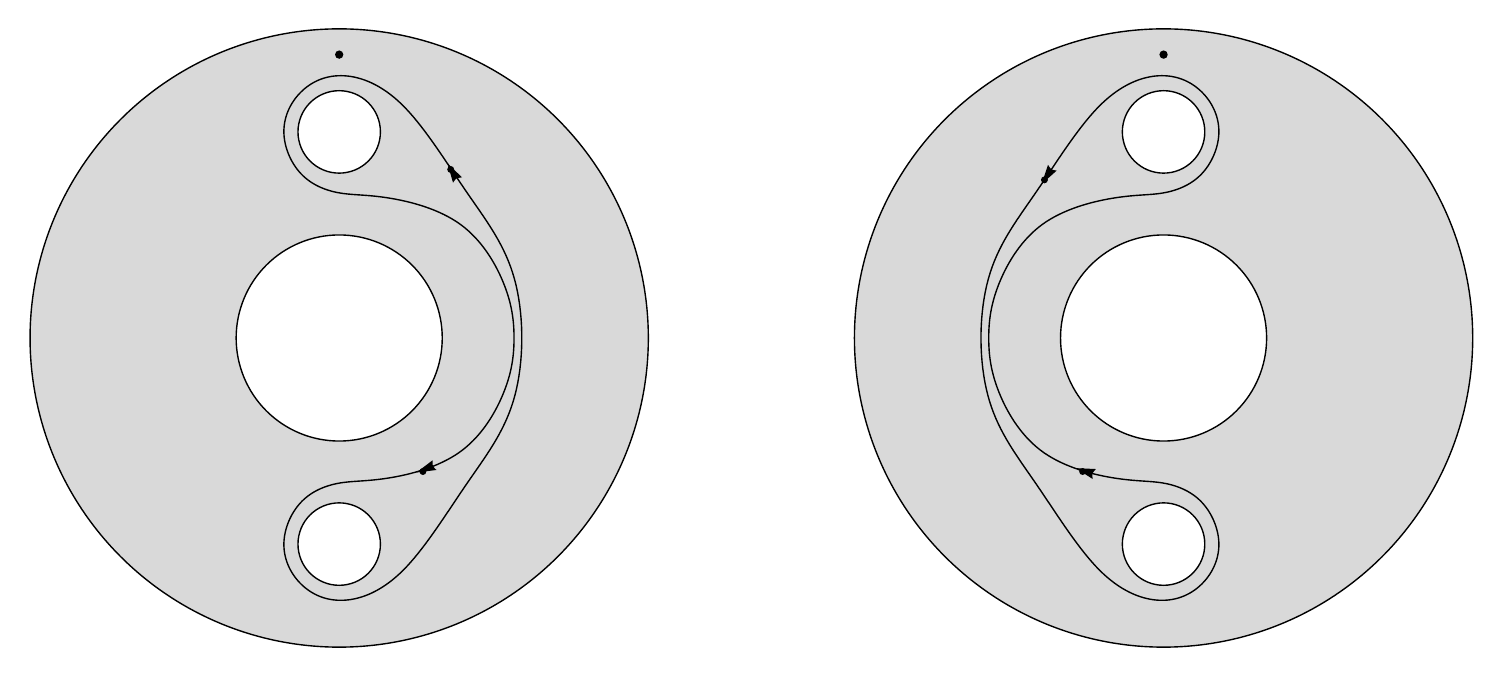}}
\unitlength1cm
\begin{picture}(0.01,0.01)
  \put(2.75,0){\footnotesize$(U,u)$  }
  
  \put(10.3,0){\footnotesize$(U,u)$  }
  \put(3.2,6.12){$\scriptstyle u$}
  \put(10.75,6.12){$\scriptstyle u$}
\end{picture}

\vspace{.4cm}

We are now in a position to define the \Car bound for a suitable family of domains. For a set $K \subset \cbar$, let us denote the spherical diameter of $K$ by $\mathrm{diam^\#}(K)$.

\begin{definition}
Let $\Uu$ be as above and (if $n \ge 2$) let $\Gamma$ be a maximally close system of meridians for $\Uu$.  We define the {\rm \Car bound} and the {\rm extended \Car bound}, $|\Uu|$ and $|\Uu|_{\mathrm E}$, respectively, of $\Uu$ by 

\vspace{-.3cm}
\begin{eqnarray*}
|\Uu| \,\,\,\,  &=& |\log ( \delta^\#(u)\, \diams (\cbar \setminus U)) | \quad \!+\! \quad  {\calL}(\Uu) \quad  \!+\!  \quad {\calD}(\Uu) ,\\
\\
|\Uu |_{\mathrm E} &=& |\log (\delta^\#(u)\, \diams (\cbar \setminus U)) | \quad \!+\! \quad \calL_{\mathrm E}(\Uu) \quad \!+\!  \quad \calD_{\mathrm E}(\Uu)
\end{eqnarray*}

\vspace{0cm}
where the terms ${\calL}(\Uu)$, $ {\calD}(\Uu)$, ${\calL}_{\mathrm E}(\Uu)$, ${\calD}_{\mathrm E}(\Uu)$ are all zero if $n=1$.

Let $n \ge 1$ and let $\calU = \{\Uaua\}_{\alpha \in A}$ be a family of pointed non-degenerate $n$-connected domains. We define the {\rm \Car bound} and the {\rm extended \Car bound} of $\calU$, $\normU$ and $\normUE$ by
\vspace{.4cm}
\[ \normU \;= sup_{\alpha} | \Uaua |, \qquad
\normUE \;= sup_{\alpha} | \Uaua |_E\]
\end{definition}

One can think of the three terms above as each preventing a
different way in which a \Car limit of non-degenerate $n$-connected domains can fail to be another non-degenerate $n$-connected domain. The first term prevents any \Car limit from simply being a point or $\C$ (recall that the spherical diameter of all of $\cbar$ is $\tfrac{\pi}{2}$ and in particular finite) and also ensures that the set of basepoints $\{u_\alpha\}_{\alpha \in A}$ is bounded. The second term prevents complementary components from merging or becoming points as in Figure 4 while the third term prevents 
components of the complement from being `engulfed' by other components, again as in Figure 4. 

Before we can state our theorem, we need more definitions relating to families of functions defined on varying domains. 

\begin{definition}
Let $\calU = \{\Uaua\}_{\alpha \in A}$ be a family of pointed hyperbolic domains in $\C$. We say $\calU$ is {\rm hyperbolically non-degenerate} or simply {\rm non-degenerate} if the limit of any convergent sequence in $\calU$ is another hyperbolic domain. 
\end{definition}

Note that a bounded family of non-degenerate pointed domains of the same connectivity is clearly hyperbolically non-degenerate. The following is immediate. 

\begin{lemma}
If $\calU = \{\Uaua\}_{\alpha \in A}$ is a non-degenerate family of pointed hyperbolic domains, then the quantities 
 $| \log(\delta^\#_{U_\alpha}(u_\alpha) \,\diams (\cbar \setminus U_\alpha ) )|$ are uniformly bounded in $\alpha$.
\end{lemma}

The converse of this is easily seen to be false, for example by considering the sequence $\{({\mathrm A}(0,\tfrac{1}{m}, m), 1)\}_{m=2}^\infty$.

\begin{definition}
If $\calF = \{f_\alpha\}_{\alpha \in A}$ is a family of analytic functions 
with each $f_\alpha$ defined on the correpsonding set $U_\alpha$ of a non-degenerate family $\calU$, then we say $\calF$ is {\rm normal} on $\calU$ if for every convergent sequence $\{(U_{\alpha_m}, u_{\alpha_m})\}_{m=1}^\infty$ in $\calU$ as above with limit $\Uu$ with $U \ne \{ u \}$, we can find a subsequence for which the corresponding functions $f_{\alpha_m}$ converge uniformly on compact subsets of $U$ (in the sense given in Definition 3.1). 
\end{definition}

\begin{definition}
Let $\calU = \{\Uaua\}_{\alpha \in A}$ be a non-degenerate family of hyperbolic pointed domains and  let $\calF = \{f_\alpha\}_{\alpha \in A}$ be a family of analytic functions 
with each $f_\alpha$ defined on $U_\alpha$. 

We say that 
$\calF$ is {\rm equicontinuous
on compact subsets of $\calU$} if for any hyperbolic distance $R >0$ there 
exists $M \ge 0$ depending on $R$ such that, for each $\alpha \in A$, $f_\alpha$ is $M$-Lipschitz with respect to the hyperbolic metric within hyperbolic distance $\le R$ of $u_\alpha$ in $U_\alpha$. In this case we write $\calF \triangleleft \,\calU$.

We say that $\calF$ is {\rm bi-equicontinuous 
on compact subsets of $\calU$} if each $f_{\alpha}$ is a (locally injective) covering map onto its image and for each hyperbolic 
radius $R >0$ there exists $K(R) \ge 1$ such that, for each $\alpha \in A$, 
$f_\alpha$ is $K$-Lipschitz and locally $K$-bi-Lipschitz 
within hyperbolic distance $\le R$ of $u_\alpha$ in 
$U_\alpha$. In this case we write $\calF \bowtie \,\calU$.
\end{definition}

Note that if all the domains in our family are the same, then equicontinuity on compact subsets and local boundedness correspond with the standard 
definitions for a family of analytic functions on some domain.

\begin{proposition}
Let $\calU = \{\Uaua\}_{\alpha \in A}$ be a non-degenerate family of hyperbolic pointed domains and  let $\calF = \{f_\alpha\}_{\alpha \in A}$ be a family of analytic functions 
with each $f_\alpha$ defined on $U_\alpha$. Then $\calF$ is normal on $\calU$ if and only if it is equicontinuous on $\calU$.
\end{proposition}

\proof If is normal on $\calU$, it follows easily that $\calF \triangleleft \,\calU$. The other direction follows easily from  ii) of \Car convergence, Montel's theorem and a standard diagonalization argument as in the proof of the Arzela-Ascoli theorem. $\Box$

\begin{definition} Let $\calU = \{\Uaua\}_{\alpha \in A}$ and $\calV = \{\Vava\}_{\alpha \in A}$ be families indexed by the same set $A$ with $\calU$ being non-degenerate. We say that a family $\calF = \{f_\alpha\}_{\alpha \in A}$ {\rm maps $\calU$ to $\calV$} if, for each $\alpha$, $f_\alpha(U_\alpha) \subset V_\alpha$ and $f_\alpha(u_\alpha) = v_\alpha$ and we write $\calF: \calU \mapsto \calV$. 

By convention, if in addition $f_\alpha$ is a covering map, we will require that $f_\alpha(U_\alpha) = V_\alpha$.
\end{definition}

Hyperbolic non-degeneracy, and local boundedness and bi-equicontinuity are related by the following result whose proof we will also postpone until after the statement of Theorem 4.2. 

\begin{theorem} Let $\calU = \{\Uaua\}_{\alpha \in A}$ be a family of pointed hyperbolic domains and let $\Pi = \{\pi_\alpha\}_{\alpha \in A}$ denote the family of normalized covering maps from $\DO$ to $\calU$. Then $\calU$ is non-degenerate if and only if $\Pi \bowtie \DO$.
\end{theorem}

Recall that by Theorem 3.1, an $n$-connected domain $U$ with $n \ge 2$ is conformally equivalent to a standard domain. This domain is unique up to rotation and to specifying which of the components correspond to the closed unit disc and the the unbounded complementary component of the standard domain. 

Such a standard domain is specified uniquely by $3n - 5$ real numbers. A pointed $n$-connected domain is then conformally equivalent to a pointed standard domain which is described using $3n - 3$ numbers. However, as the pointed standard domain for an $n$-connected domain is unique only up to rotation, this allows us to eliminate one more parameter so that it is specified uniquely using $3n - 4$ real numbers (where, for example, we insist that the basepoint of the standard domain lie on the positive real axis).

We consider the set $\Dn$ of all ordered 4-tuples $(U, u, K^1, K^2)$, where $U$ is a non-degenerate $n$-connected  domain, $u \in U$ and $K^1$, $K^2$ are two distinct components of the complement $\cbar \setminus U$. 
A sequence $\{(U_m,  u_m, K^1_m, K^2_m)\}_{m=1}^\infty$ is said to converge to a limit $(U,  u, K^1, K^2) \in \Dn$ if the pointed domains $(U_m,u_m)$ converge in the \Car topology to $\Uu$, and any limit in the Hausdorff topology of the sets $K^1_m$ is contained in $K^1$ and also any limit in the Hausdorff topology of the sets $K^2_m$ is contained in $K^2$. We can then use this notion of convergence to endow $\Dn$ with a corresponding topology. Conformal equivalence in $\Dn$ is an equivalence of pointed domains which also preserves the labelling of boundary components, which by Theorem 3.1 is a well-defined equivalence relation. 

The quotient of $\Dn$ by this relation gives rise to a moduli space $\Mn$, which we then endow with the quotient topology arising from the topology defined on $\Dn$. As we have seen, non-degenerate $n$-connected pointed slit domains may be identified with a subset $\Bn$ of $\R^{3n - 4}$ which is easily seen to be open and which we endow with the topology of $\R^{3n - 4}$. 
The canonical injection takes points in $\Bn$ to pointed slit domains and if we choose the first labelled complementary component of each slit domain to be $\overline \D$ and the second the unbounded complementary component, then this gives us an injection from $\Bn$ to $\Dn$ which is easily seen to be continuous. These points in $\Dn$ are then in turn mapped to points in the moduli space $\Mn$ by the above equivalence relation. 

Denoting this composition by $\inn$, it follows from the uniqueness part of Theorem 3.1 that $\inn$ is injective. Also, $\inn$ is surjective from the existence part of Theorem 3.1.  Using conformal mappings to standard domains which preserves the above labelling of complementary components shows that sets of equivalence classes in $\Mn$ for which the numbers associated with the corresponding standard domains lie in a closed subset of $\Bn$ are themselves closed (note that the notion of labelling and convergence in $\Dn$ is exactly the same as that used at the start of the proof of this theorem). From this it follows easily that 
$\Mn$ is Hausdorff as a topological space.  

As $\inn$ is a continuous bijection from $\Bn$ to the Hausdorff space $\Mn$, it is then a homeomorphism on compact subsets of $\Bn$. This has the benefit that we can identify the topology of $\Bn$ on compact subsets with that of compact subsets of $\Mn$ and thus consider $\Mn$ as a $(3n-4)$-dimensional manifold (the author gratefully  
acknowledges the help of Adam Epstein, Vaibhav Gadre and Jeremy Kahn with the above discussion). This allows us to make the following definition. 

\begin{definition} Let $\calU = \{\Uaua\}_{\alpha \in A}$ be a family of $n$-connected non-degenerate pointed domains. We say $\calU$ is {\rm bounded in moduli space} if for each domain $U_\alpha$ and any choice of distinct components $K^1_\alpha$, $K^2_\alpha$ of $\cbar \setminus U_\alpha$ which gives us an embedding into $\Dn$, the image of this family in $\Mn$ forms a set whose inverse image under $\inn$ is precompact in $\Bn$.
\end{definition}

We remark that the subset of $\Mn$ which arises for such a family must also be precompact. Before stating the main result of this paper, we remind the reader that one condition is said to imply another up to constants if the constants for the first condition imply non-trivial bounds on those for the second. 

\begin{theorem}Let $\calU = \{\Uaua\}_{\alpha \in A}$ be a family of pointed $n$-connected\\ domains. If $n \ge 2$, then the following are equivalent:
\begin{enumerate}
\item $\pt \sqsubset \calU \sqsubset \cbar$;

\vspace{.25cm}
\item $\normUE \;< \infty$;

\vspace{.25cm}
\item $\normU \;< \infty$;

\vspace{.25cm}
\item There exist constants $\delta_1, \delta_2 > 0$ and curves $\eia$, $ \alpha \in A, 1 \le i \le P(n)$ for which we have the following:

\vspace{.2cm}

\renewcommand{\theenumii}{\emph {\alph{enumii})}}
\renewcommand{\labelenumii}{\theenumii}

\begin{enumerate}

\item $\eta^i_\alpha$ is a curve in $U_\alpha$ and $u_\alpha$ lies on each $\eta^i_\alpha$;

\vspace{.2cm}
\item $\eia$ separates $K^i_\alpha$ from the other components of $\cbar \setminus U_\alpha$;

\vspace{.2cm}
\item Any point on $\eia$ is at least spherical distance $\delta_1$ away from $\cbar \setminus U_\alpha$;

\vspace{.2cm}
\item The spherical diameters of the components $K^i_\alpha$, $1 \le i \le n$, of $\cbar \setminus U_\alpha$ are at least $\delta_2$;

\end{enumerate}

\vspace{.25cm}
\item If for each $\alpha$ we let $(A^{\Lambda_\alpha}, a_\alpha)$ be a pointed standard domain for $(U_\alpha, u_\alpha)$ (where for $n \ge 2$ we allow any choice of which components of $\cbar \setminus U_\alpha$ correspond to the closed unit disc and unbounded complementary component of $A^{\Lambda_\alpha}$), then the corresponding family $\calA = \{(A^{\Lambda_\alpha}, a_\alpha)\}_{\alpha \in A}$ satisfies $\normA < \infty$ and the family of inverse Riemann maps $\psi_\alpha$ gives a univalent family $\Psi$ with $\Psi : \calA \mapsto \calU$ and $\Psi \bowtie \calA$; 

\vspace{.25cm}
\item $\calU$ is hyperbolically non-degenerate and bounded in moduli space. 

\end{enumerate}
Furthermore, conditions 2., 3., 4. and 5. are equivalent up to constants.

If $n=1$, part b) of condition 4. is vacuous while a) and c) become the single condition $\delta^\#_{U\alpha}(u_\alpha) \ge \delta_1$. For 5., the standard pointed domains are all $\DO$. 
\end{theorem}

{\bf Proof of Lemma 4.1\hspace{.4cm}} The conclusion is trivial for $n = 2$ as the equator of a conformal annulus is unique and so we can assume $n \ge 3$. Again for convenience, let us assume that $U \subset \C$ so that we can use homology to characterize how simple closed curves separate $\cbar \setminus U$. 

Let $2 \le i \le n-2$, select integers $1 \le j_1 < j_2 < \cdots\cdots < j_i \le n$  and let  $K^{j_1}, K^{j_2}, \ldots\ldots, K^{j_i}$ be (distinct) components of $\cbar \setminus U$. Selecting components in this way is equivalent to choosing a homology class of a simple closed curve in $U$ where we separate these components of the complement of $U$ from the remaining ones. Note that we do not need to consider the cases where $i =1$ or $n-1$ as these are the cases where we have principal meridians in which case the result is trivially true (this includes the degenerate case where some of the principal meridians fail to exist). We then set 
\vspace{.2cm}
\[d = \inf_\gamma \rho(u, \gamma)\]

where the infimum is taken over all those meridians which separate the components $K^{j_1}, K^{j_2}, \ldots, K^{j_i}$ from the rest of $\cbar \setminus U$.  Let $\{\gamma_m\}$ be a sequence of meridians which are in this homology class and for which the distances $d_m = \rho_U( u, \gamma_m)$ converge to $d$. 

Now let $\pi: \D \to U$ with $\pi(u) = 0$, $\pi'(u) > 0$ be the normalized covering map as in Theorem 1.2 and for each $m$ let $\sigma_m$ be a lift of $\gamma_m$ which is at hyperbolic distance $d_m$ from $0$ in $\D$. Next let $\eta_m$ be a full segment of $\gamma_m$ of length $\ell_m = \ell_U(\gamma_m)$ which is at distance $d_m$ from $0$ so that $\pi(\eta_m) = \gamma_m$ and $\pi$ is injective on $\eta_m$ except at the endpoints. Note that by Theorem 1.6 again, all the segments $\eta_m$ will have the same hyperbolic length. 

Passing to a subsequence if necessary, we can choose our segments $\eta_m$ so that 
they converge to a segment $\eta$ of another geodesic $\sigma$ which 
is at distance $d$ from $0$. If we now let $\gamma = \pi(\eta)$, then it is clear 
that $\gamma$ is a closed geodesic in $U$ and similarly to in the proof of Theorem 1.8, one sees that $\gamma$ is a smooth curve without any `corners'. 

It also follows from the uniform convergence of the segments $\eta_m$ to $\eta$ that $\gamma_m$ is homotopic to $\gamma$ for $m$ large enough. By Theorem 7.2.5 on Page 129 of \cite{KL}, $\gamma_m = \gamma$ for $m$ large enough and it follows that $\gamma$ must be simple and a meridian in the desired homology class. Repeating this argument for each choice of subset of components of $\cbar \setminus U$ then gives the required extended system of meridians. $\Box$

{\bf Proof of Theorem 4.1\hspace{.4cm}} The direction where $\calU$ is non-degenerate follows by Theorem 1.2.
Suppose now that $\Pi \bowtie \DO$, let $\{U_{\alpha_m}, u_{\alpha_m}\}_{m=1}^\infty$ be a sequence which converges in the \Car topology to a limit $\Uu$ and suppose for the sake of contradiction that $U$ is not a hyperbolic domain. As $\Pi \bowtie \DO$ by Proposition 4.1, $\Pi$ is normal and it follows from Bloch's Theorem (e.g. \cite{Con} Page 293, Chapter XII, Theorem 1.4) that $U$ cannot be a point. By postcomposing as usual with a suitable M\"obius transformation, we can assume without loss of generality that $U$ is either $\C \setminus \{0\}$ or all of $\C$.

If $U = \C \setminus \{0\}$, on passing to a subsequence if needed, we can clearly assume that the covering maps $\pi_{\alpha_m}$ converge to a limit function $\pi$. By the bi-equicontuity of $\Pi$, $\pi$ must be locally injective. Using Rouch\'e's theorem and local compactness as observed by Epstein in the remarks preceding the proof of Lemma 6 on Page 15 of \cite{Ep}, it follows that $\pi(\D) = \C \setminus \{0\}$. 

Now $e^z$ maps $\C$ to $\C \setminus \{0\}$ and, by the monodromy theorem, we can then find a branch of $\pi^{\circ -1}( e^z)$ which is defined on all of $\C$. However, this gives us an entire function whose range is a subset of $\D$ which contradicts Picard's theorem on the range of an entire function (e.g. page 319 of Lang's book \cite{Lang}).  Finally, the argument where $U = \C$ is then a simpler version of this.  $\Box$

{\bf Proof of Theorem 4.2\hspace{.4cm}}For \emph{1.} $\Longrightarrow$ \emph{2.}, if the first term in the \Car bound for $\calU$ was unbounded, using the Hausdorff version of \Car convergence, we could find a sequence which tended to either a point (which could possibly be infinity as we defined \Car convergence using convergence of the base points in the spherical and not the Euclidean topology) or all of $\cbar$, both of which would contradict the boundedness of $\calU$. The required bounds on the lengths and distances of the pointed domains of $\calU$ follow again from the boundedness of $\calU$ and Corollary 1.1. 

\emph{2.} $\Longrightarrow$ \emph{3.} is immediate. 

For \emph{3.} $\Longrightarrow$ \emph{4.}, for each $\alpha$ and each $1 \le i \le P(n)$, let $\eia$ be the curve obtained by adding to the principal meridian $\gia$ the hyperbolic segment which connects $u_\alpha$ to the closest point on $\gia$ in the hyperbolic metric and which is then traversed in both directions (note that we can easily modify $\eia$ slightly if we wish to ensure that it is a simple closed curve). $u_\alpha$ then lies on $\eia$ and it is clear that this curve separates $K^i_\alpha$ from the other components of $\cbar \setminus U_\alpha$ and so we have \emph{a)} and \emph{b)}.

We will need to apply the estimates of Lemma 3.1 on the hyperbolic metric in a uniform manner which will be independent of the choice of pointed domain in $\calU$ and of the point at which we wish to estimate the hyperbolic metric for that domain. In particular we want to find $C > 0$ and $\delta > 0$ such that for any $\alpha$ and any $z \in \partial U_\alpha$, these estimates hold for this $C$ for all $w \in U_\alpha$ within spherical distance $\delta$ of $z$. 

In view of our remarks after the statement of Lemma 3.1, we will be able to do this if we can prove the claim that we can find $d >0$ independent of $\alpha$ such that if $z^1_\alpha$ is any point in $\cbar \setminus U_\alpha$, we can find two other points $z^2_\alpha$ and $z^3_\alpha$ also in $\cbar \setminus U_\alpha$ so that the these three points are separated by at least distance $d$ from each other in the spherical metric. 

So suppose this fails. This means that either we can find a sequence in $\calU$ which tends $\cbar$ with either one or two points removed. As usual, without loss of generality, we can assume such a sequence tends either to $\C$ 
or to a punctured plane in which case the corresponding domains contain an annulus which separates the components of the complement and whose modulus tends to infinity. 

The first case is excluded in view of the first term of the \Car bound. 
In the second case, 
the equator of such an annulus has very small hyperbolic length and by Theorem 1.6 separates the complement of the domain in the same way as a meridian whose hyperbolic length is at least as small. If the connectivity $n$ of all the domains in $\calU$ is $2$ or $3$, then this is impossible in view of the second term in the \Car bound as here all meridians are principal. 

On the other hand, if all the domains of $\calU$ have higher connectivity and thus have meridians which are not principal, it follows from Theorem 2.6 of \cite{Com2} that the principal meridians do not intersect any other meridians of a domain (including other principal meridians). Additionally, a meridian which is not principal must contain principal meridians in both of its complementary components. By the collar lemma (\cite{KL} Page 148, Lemma 7.7.1), we can then find a principal meridian whose distance to this very short meridian is unbounded and this is impossible in view of the third term in the \Car bound. The claim then follows.  

The first term of the \Car bound implies that the distances $\delta^\#(u_\alpha)$ to the boundary are uniformly bounded below and the second and third terms together imply that the hyperbolic lengths of the curves $\eia$ are uniformly bounded above. As we have shown the lower estimate for the hyperbolic metric from Lemma 3.1 is uniform, and as the 
improper integral
\vspace{.16cm}
\[ \int \limits_0^{\:\:\:\tfrac{1}{2}} { \frac{1}{x \log (1/ x)}\,{\rm d}x}\]

diverges, we can deduce that we can find $\delta_1$ not depending on $\alpha$ such that the curves $\eia$ are at least (spherical) distance $\delta_1$ away from the complement of $U_\alpha$ as desired and so we have shown \emph{c)}. 

Finally to prove {\emph d)}, we first note that from above the complementary components are all distance at least $2\delta_1$ apart. Thus if one of the components $\Kia$ of the set $\cbar \setminus U_\alpha$ had arbitrarily small spherical diameter, then $U_\alpha$ would contain a (round) annulus of very large modulus separating $\Kia$ from the rest of $\cbar \setminus U_\alpha$. If we then consider the principal meridian $\gamma^j_\alpha$ (where $j=1$ if $n = 2$ and $j = i$ if $n > 2$), then $\gamma^j_\alpha$ separates $\Kia$ from the rest of $\cbar \setminus U_\alpha$ as does the equator of this annulus. 

By Theorem 1.6, the length $l^j_\alpha$ of $\gamma^j_\alpha$ would then be very small. However, the third term in the \Car bound shows that the lengths $l^j_\alpha$ are uniformly bounded below and so we must have $\delta_2 >0$ independent of $\alpha$ for which all the complements $\Kia$ of $\cbar \setminus U_\alpha$ have (spherical) diameter at least $\delta_2$ as required. 

To show \emph {4.} $\Longrightarrow$ \emph{1.}, suppose we have $\delta_1, \delta_2 > 0$ and curves $\eia$ as above and suppose we have a sequence which converges in the 
\Car topology which we will label $\{\Umum\}_{m=1}^\infty$. Let the limit of this sequence be $\Uu$. By conditions \emph{a)}  and \emph {c)} of \emph {4.} $\delta^\#_{U_\alpha}(u_\alpha) \ge \delta_1$ for every $\alpha$. It then follows easily using iii) of \Car convergence that $U \ne \{u\}$. 

Now take a subsequence $m_k$ so that  complementary components of each $\Umk$ converge in the Hausdorff topology. By \emph{4.b),c)} the Hausdorff limit must consist of exactly $n$ components $K^1, K^2, \ldots \ldots, K^n$ which are at least (spherical) distance $2\delta_1$ apart and by \emph{4.d)} each of these must have diameter at least $\delta_2$.

We need to show that each of the sets $K^i$ lies in a different component of $\cbar \setminus U$. By relabelling if necessary,  we can ensure that for each $i$  the sets $\Kimk$ converge in the Hausdorff topology to $K^i$. As a \Car limit of $n$-connected domains, $U$ has connectivity $n' \le n$ and we need to show $U$ has connectivity exactly $n$ and that each component of $\cbar \setminus U$ contains precisely one of the sets $K^i$ above. 

We first show that for every $i$ we can find a point $\partial K^i$ which is in $\partial U$. Fix $1 \le i \le n$. The curve $\eta^i_{m_k}$ separates $K^i_{m_k}$ from the rest of $\cbar \setminus U_{m_k}$. By joining the closest points on this curve and $K^i_{m_k}$ with a line segment, we can extend this to a curve $\tilde \eta^i_{m_k}$ which joins $u_{m_k}$ with a point on $\partial K^i_{m_k}$. 

Now pass to a further subsequence if needed so that these curves converge in the Hausdorff topology to a continuum $\tilde \eta^i$ which by i) of \Car convergence then clearly joins $u$ to $K^i$. By making this curve a little smaller if needed, we can assume that it meets $K^i$ only at its endpoint and so terminates in a line segment of (spherical) length $\ge \delta_1$ at some point $\tilde z \in \partial K^i$. 

From iii) of \Car convergence it then follows that any compact subset of  $\tilde \eta^i \setminus \{\tilde z\}$ will lie in $U$ whence $\tilde \eta^i \setminus \{\tilde z\} \subset U$. $\tilde z$ can then be approximated by points in $U$ and since $\tilde z \in K^i$, it follows from the Hausdorff convergence of the sets $\Kimk$ to $K^i$ and ii) of \Car convergence that $\tilde z$ cannot be in $U$. Thus $\tilde z \in \partial K^i \cap \partial U$ as we want.




By Lemma 2.1 if $z \in \partial U$, then $z$ must meet $K^i$ for some $i$ and since the sets $K^i$ are at least distance $2\delta_1$ apart, there can be only one such $i$. Now if $U$ had connectivity $ < n$, since from above each of the sets $\partial K^i$ meets $\partial U$, we could find a component $L$ of $\cbar \setminus U$ and two points $z_1$, $z_2$ of $\partial L \subset \partial U$, which were contained in sets $K^{i_1}$, $K^{i_2}$ respectively with $i_1 \ne i_2$. Now for each $1 \le i \le n$, let $G^i = \partial L \cap K^i$. The sets $G^i$ are each clearly closed and from above they give a non-trivial separation of $\partial L$. However, by the corollary to Theorem 14.4 on Page 124 of \cite{New}, $\partial L$ is connected and with this contradiction we see that $U$ must be $n$-connected. 

Now for each $1 \le i \le n$, let $L^i$ be the component of $\cbar \setminus U$ whose boundary meets $K^i$. As $K^i$ is a connected subset of $\cbar \setminus U$, we must have that $K^i \subset L^i$. It then follows from above that the sets $L^i$ will have spherical diameter $\ge \delta_1$. Thus the \Car limit $\Uu$ is a non-degenerate $n$-connected pointed domain and so $\calU$ is bounded. We still need to show that \emph{4.} $\Longrightarrow$ \emph{2.} holds up to constants. However, the bound on the first term of the extended \Car bounds follows directly from what we have just proved while the bound on the other two terms follows from what we have just proved combined with Corollary 1.1. 

We now show \emph{2.} $\Longleftrightarrow$ \emph{5.} up to constants. So suppose \emph{2.} holds. It follows from \emph{2.} $\Longrightarrow$ \emph{1.} and Theorem 3.2 that 
for the points $a_\alpha = \phi_\alpha (u_\alpha)$, 
the quantity $\log (\delta^\#_{A^{\Lambda_\alpha}}(a_\alpha))$ is uniformly bounded in $\alpha$ and that this bound is uniform with respect to $||\,\calU ||$. Note that, for $n \ge 2$, on examining the start of the proof of Theorem 3.2, we see that by passing to a subsequence if necessary, we can label the complementary components of a limit pointed domain and the corresponding choice of standard domain to be compatible with those of the approximating pointed domains as specified in the statement of Theorem 3.2. Now there are only $n(n-1)$ ways of assigning which complementary components of a limit domain are mapped to the unit disc or the unbounded complementary component of a standard domain. Hence by Theorem 3.2 and ii) of \Car convergence applied to the standard domains, the bound on $\log (\delta^\#_{A^{\Lambda_\alpha}}(a_\alpha))$ will not depend on our choice of standard domain for each $U_\alpha$.

Since for $n \ge 2$ all standard domains avoid $\overline \D$, and for $n=1$, they avoid $\cbar \setminus \D$, the quantity $\log (\diams (\cbar \setminus A^{\Lambda_\alpha}))$ is also uniformly bounded in $\alpha$ and must also be uniform with respect to $||\,\calU ||$. Thus the first term of the \Car bound is uniformly bounded for the family $\calA$ and this bound is uniform with respect to $||\,\calU ||$. By Lemma 3.2, meridians and principal meridians are conformally invariant and using this and the 
conformal invariance of hyperbolic length, the second and third terms for each $|\ALaa|$ are the same as those for $|\Uaua|$ and with this we have shown $||\hspace{.025cm}\calA\hspace{.025cm}|| < \infty$ and that this bound is uniform with respect to $||\,\calU \hspace{.025cm}||$. 

We still need to show that $\Psi \bowtie \calA$ and that the estimates on $\Psi$ are uniform with respect to $||\,\calU\hspace{.025cm}||$. For each $\alpha$ let $\chi_\alpha$ be the unique normalized covering map from $\DO$ to $\ALaa$ as in Theorem 1.2. As in the proof of Theorem 3.2, the composition $\pi_\alpha = \psi_\alpha \circ \chi_\alpha$ is then the corresponding normalized covering map for $\Uaua$.  

Now let $\{(U_{\alpha_m},u_{\alpha_m})\}_{m=1}^\infty$ be any sequence in $\calU$ with corresponding pointed standard domains $\{(A^{\Lambda_{\alpha_m}}, a_{\alpha_m})\}_{m=1}^\infty$. Using \emph{2.} $\Longrightarrow$ \emph{1.} once more, by selecting a subsequence if necessary, we can ensure that both sequences converge to non-degenerate $n$-connected pointed domains $\Uu$ and $(A^\Lambda,a)$ respectively and where of course $A^\Lambda$ must be a standard domain. 

By Theorem 1.2, if we let $\pi$ and $\chi$ be the normalized covering maps for $U$ and $A^\Lambda$ respectively, then $\pi_m \to \pi$ and $\chi_m \to \chi$ locally uniformly on $\D$. It then follows that for each hyperbolic radius $R \ge 0$, we can find $K = K(R) >1$ not depending on $\alpha$ such that within hyperbolic distance $R$ of $a_\alpha$ in $A^{\Lambda_\alpha}$, $\psi_\alpha$ is locally bi-Lipschitz 
with constant $K$. The same argument shows that $K$ is uniform with respect to $R$ and $||\,\calU\hspace{.025cm}||$ and so we have shown \emph{2.} $\Longrightarrow$ \emph{5.} up to constants.  

To show \emph{5.}$\Longrightarrow$ \emph{2.} up to constants, as $||\hspace{.025cm}\calA\hspace{.025cm}|| < \infty$ and $\Psi \bowtie \, \calA$, the quantities $\log \delta^\#_{U_\alpha}(u_\alpha)$ are obviously uniformly bounded below on applying the Koebe one-quarter theorem. Again for the same reasons as before, the second and third terms for each $|\Uaua|$ are the same as those for $|\Aaaa|$. 

Suppose now we could find a sequence $\{(U_{\alpha_m}, u_{\alpha_m})\}_{m=1}^\infty$ with limit $(\cbar \setminus \{v\}, u)$ for some points $u, v \in \cbar$ with $u \ne v$. As usual, for convenience we can assume that $v = \infty$ so that our sequence converges to $(\C, u)$. Now let $\phi_m$ be the Riemann map from each $U_m$ to the corresponding standard domain. By ii) of \Car convergence, \emph{2.} $\Longrightarrow$ \emph{1.} applied to $||\hspace{.025cm}\calA\hspace{.025cm}|| < \infty$ and Montel's theorem, this would clearly 
give a normal family on any bounded open subset of $\C$. Any limit function would then be an entire function which was either constant or a univalent entire function which avoided either $\overline \D$ (for $n \ge 2$) or $\C \setminus \D$ (for $n = 1$). The case of a constant limit function would violate the bi-equicontinuity of $\Psi$ while the non-constant case would violate Picard's theorem on the range of a non-constant entire function. Thus the quantity $\log(\diams(\cbar \setminus U_\alpha))$ is uniformly bounded and this bound is uniform with respect to the bounds for $\calA$ and the bi-equicontinuity of $\Psi$.

Finally, we show \emph{5.} $\Longleftrightarrow$ \emph{6.} Suppose first \emph{5.} holds. From \emph{5.} $\Longrightarrow$ \emph{2.}, we have that $\calU$ is hyperbolically non-degenerate, while the fact that $||\hspace{.025cm}\calA\hspace{.025cm}|| < \infty$ and  \emph{2.} $\Longrightarrow$ \emph{1.} for $\calA$ shows that 
$\calU$ must be bounded in moduli space. For the other direction, if \emph{6.} holds, then by Lemma 4.2, as $\calU$ is non-degenerate, the first term in the \Car bound for the pointed domains of $\calU$ is uniformly bounded. Since $\calU$ is bounded in moduli space, if we label two of the complementary components of each $U_\alpha$ in any way we like, 
the vectors associated with the correpsonding family of standard domains form a precompact set in $\Bn$ whence the corresponding standard domains are precompact in $\Dn$. This shows that the resulting family $\calA$ of standard domains is bounded. The uniform bounds on the other two terms of the \Car bounds for the domains of $\calU$ then follows from Lemma 3.2 on the conformal invariance of principal meridians together with \emph{2.} $\Longrightarrow$ \emph{5.} for $\calU$. With this the proof is complete. $\Box$

Theorem 4.2 has a useful consequence.

\begin{corollary}
Let $\calU = \{\Uaua\}_{\alpha \in A}$ be a hyperbolically non-degenerate family of non-degenerate pointed domains each of which has finite connectivity such that the 
the the hyperbolic length and distance to the base point for all the meridians of each $U_\alpha$ are uniformly bounded above. Then there exists $K_1 > 0$ independent of $\alpha$ for which the boundaries $\partial U_\alpha$ are $K_1$-uniformly perfect. Equivalently, there exists $K_2 \ge 1$ such that if $\rho_\alpha (\cdot \,,\cdot)$ denotes the hyperbolic metric on $U_\alpha$, then
\vspace{.25cm}
\[ \frac{1}{K_2} \frac{\dsharpz}{\delta^\#_{U_\alpha}(z)} \le {\rm d}\rho_\alpha(z) \le K_2 \frac{\dsharpz}{\delta^\#_{U_\alpha}(z)}.\]

In particular, the above holds if $\calU$ is a bounded family.

\end{corollary}

\proof  The uniform bound on the lengths and distances of each $U_\alpha$ implies that the lengths of all the meridians of each $U_\alpha$ (and not merely the principal meridians) are uniformly bounded. To see this, note that by the collar lemma, if we could find arbitrarily short meridians, then the same argument as used in \emph{3.} $\Longrightarrow$ \emph{4.} of the last result would imply that the distances for $\calU$ were unbounded. The bound on the moduli of annuli which separate the complements of the domains 
follows from the upper bound on the hyperbolic lengths of all the meridians, Theorem 1.6 and the fact that conformal annuli of large modulus contain genuine annuli of large modulus. 

The above condition on the hyperbolic metric is usually given in terms of the Euclidean distance to the boundary and follows from \cite{BP} Corollary 1. To see that we also have the above version using the spherical metric, we observe that 
it follows from the hyperbolic non-degeneracy, the condition on the lengths and a discussion similar to that in the proof of {\emph 3.} $\Longrightarrow$ {\emph 4.} in Theorem 4.2 about applying the estimates of Lemma 3.1 in a uniform manner that we can find three points in the complement of each domain whose separation in the spherical metric is uniformly bounded below. Theorem 2.3.3 on page 34 of \cite{Bear} shows that the M\"obius transformations which map these three points to $0$, $1$ and $\infty$ then give a bi-equicontinuous family as defined in Section 1 of the first part of this paper and the desired estimates then follow easily in the same manner as those for Lemma 3.1. $\Box$

Note that it is essential we consider either the lengths and distances or just the extended lengths in the above statement as the following example shows. For each $m \ge 3$, let $U_m$  be the quadruply connected domain obtained by removing the two small closed discs $\overline {\mathrm D}(\tfrac{5}{4m}, \tfrac{3}{4m})$, $\overline {\mathrm D}(-\tfrac{5}{4m}, -\tfrac{3}{4m})$ and two large closed discs $\overline {\mathrm D}(\tfrac{5m}{4}, \tfrac{3m}{4})$, $\overline {\mathrm D}(-\tfrac{5m}{4}, -\tfrac{3m}{4})$ from $\cbar$ and consider the sequence of pointed domains $\{(U_m, 1)\}_{m=3}^\infty$ (note that these discs are chosen so that each $U_m$ is invariant under the transformation $z \mapsto \tfrac{1}{z}$). 

By rescaling in turn about $0$ and $\infty$, we can apply the estimates of Lemma 3.1 in a uniform manner to deduce that the lengths of the principal meridians are uniformly bounded above and below away from $0$ in $m$. If we then let $\gamma_m$ be the meridian which separates $\overline {\mathrm D}(\tfrac{5}{4m}, \tfrac{3}{4m}) \cup \overline {\mathrm D}(-\tfrac{5}{4m}, -\tfrac{3}{4m})$ from $\overline {\mathrm D}(\tfrac{5m}{4}, \tfrac{3m}{4}) \cup \overline {\mathrm D}(-\tfrac{5m}{4}, -\tfrac{3m}{4})$, then it follows by comparing the length of this meridian with that of the unit circle using Theorem 1.6 that the 
length of this meridian clearly tends to $0$ as $m$ tends to infinity. The point here is that the lack of uniform bounds on the moduli of annuli which separate the complements of these domains is seen only in terms either of
the fact that the two pairs of principal meridians for these domains are getting farther apart or of the behaviour of a meridian which is not a principal meridian for any of the pointed domains in this sequence. 

\section{Families of Functions}

In this section we present some theorems which illustrate how the results we have proved show how the the definitions we gave of normality, equicontinuity and local uniform convergence are natural extensions of the standard ones. We also present some further useful results concerning bounded families of pointed domains.

From now on let us assume that $\calU = \{\Uaua\}_{\alpha \in A}$ and $\calV = \{\Vava\}_{\alpha \in A}$ are two families of pointed domains indexed by the same set $A$, and let $\calA = \{\ALaa\}_{\alpha \in A}$ denote a family of standard pointed domains for $\calU$. 

If we let $\calF =  \{f_\alpha\}_{\alpha \in A}$, $\calG =  \{g_\alpha\}_{\alpha \in A}$ be families defined on $\calU$, $\calV$ respectively with $\calF : \calU \mapsto \calV$, then we obtain the family of compositions 
$ \{g_\alpha \circ f_\alpha \}_{\alpha \in A}$ which is a family defined on $\calF$ and which we denote in the obvious way by $\calG \circ \calF$. Similarly, if the members of $\calF$ are all univalent, then the family of inverse functions will be defined on $\calV$ and we will denote it by $\calF^{\circ -1}$.

The following is immediate in view of the Schwarz lemma for the hyperbolic metric. 

\vspace{.2cm}
\begin{proposition}Let $\calU$, $\calV$ be two normal families and let $\calF$ be a family of functions which maps $\calU$ to $\calV$.
\vspace{0cm}

\begin{enumerate} 

\item If $\calF \triangleleft \, \calU$ and $\calG \triangleleft \, \calV$, then $\calG \circ \calF \triangleleft  \, \calU$. 

\vspace{.2cm}
\item If $\calF \bowtie \,\calU$ and $\calG \bowtie \calV$, then $\calG \circ \calF \bowtie \,\calU$. 

\vspace{.2cm}
\item If $\calF$ is univalent and $\calF \bowtie \,\calU$, then $\calF^{\circ -1} \bowtie \calV$

\end{enumerate}

\end{proposition}

Let $\Pi = \{\pi_\alpha\}_{\alpha \in A}$ denote the family of normalized covering maps from $\DO$ to $\calU$. The following is an immediate consequence of Theorems 4.1 and 4.2.  
 
\begin{proposition} If $\pt \sqsubset \calU \sqsubset \cbar$, then $\Pi \bowtie \DO$. 
\end{proposition}

The converse of this result is false. To see this consider the sequence $\{(\Um,0)\}_{m=4}^\infty$ where $U_m = \D \setminus \overline {\mathrm D}(1-2/m, 1/m)$. By Theorem 1.2 again, the normalized covering maps clearly converge to the identity and so give a bi-equicontinuous family on $\DO$. However, the limit of the doubly connected domains $U_m$ is the unit disc which is simply connected and so this sequence is not bounded. Nevertheless, for Riemann maps  and simply connected domains, the converse is true as can be seen immediately from Theorem 4.2.

From Theorem 4.1 and Proposition 5.1, we immediately get the following. 

\begin{corollary}[Bounded Families and Equicontinuity]\hspace{1cm}\\
Let $\,\calU$ and $\calV$ be families of pointed hyperbolic domains with $\calU$ hyperbolically non-degenerate and let $\calF: \calU \mapsto \calV$ be a family of covering maps. Then $\calV$ is hyperbolically non-degenerate if and only if $\calF \bowtie \, \calU$.
\end{corollary}

In view of the counterexample above, even if $\calU$ is bounded and $\calF: \calU \mapsto \calV$ with $\calF \bowtie \,\calU$, we cannot say that $\calV$ is bounded. However, we do have the following.

\begin{theorem}[Bi-Equicontinuity and Boundedness]\hspace{1cm}\\
Let $d \ge 1$, let $\calU$ be bounded and let $\calF : \calU \mapsto \calV$ with $\calF \bowtie \, \calU$ such that the degrees of the mappings of $\calF$ are uniformly bounded above by $d \ge 1$ and all the domains in $\calV$ have the same connectivity. Then $\calV$ is also bounded. Further $\normV$ is uniformly bounded with respect to $\normU$ and $d$.
\end{theorem}

In view of what we said above, the assumption on the boundedness of the degrees is certainly necessary and we proceed by first proving the following technical lemma. The main reason we require this lemma is that, although a lifting of a meridian under a covering map of finite degree must be a simple closed geodesic, there is no guarantee that it will still be a meridian. For two subsets $E,F$ of $\C$, let us use ${\dist}(E,F)$ to denote the Euclidean distance from $E$ to $F$.

\begin{lemma}Let $U$ be a hyperbolic domain, suppose $\cbar \setminus U = E \cup F$ where $E$ and $F$ are closed disjoint non-empty subsets of $\cbar$ with $\infty \in F$
and let $\tilde \gamma$ be a curve in $U$ which separates $E$ and $F$.
Let $r, R$ be two positive real numbers such that ${\dist}(\tilde \gamma, \partial U) \ge r$, $\ell(\tilde \gamma) \le R$ and $\tilde \gamma \subset \overline {{\mathrm D}(0,R)}$. Then we can find a curve $\gamma$ in $U$ with $n(\gamma, z) = 1$ for all $z \in E$, $n(\gamma, z) = 0$ for all $z \in F$ and whose hyperbolic length $\ell(\gamma)$ is uniformly bounded above in $r$ and $R$.
\end{lemma}

\proof  We begin by defining an enlarged version $\tilde E$ of $E$ by setting $\tilde E = \{z : \dist(z, E) \le \tfrac{r}{4}\}$. Now   
cover the plane with a mesh of closed squares of side length $\tfrac{r}{4\sqrt 2}$ whose boundaries are oriented positively. Let $\Gamma$ be the cycle
\vspace{.2cm}
\[\Gamma = \sum_i Q_i\]

where the sum ranges over all those squares $Q_i$ which meet $\tilde E$. Note that the number of such squares is uniformly bounded above in $r$ and $R$, a crude bound being 
$\left (2\tfrac{R}{r} + 2 \right )^2$.

By cancelling those segments which are on the boundary of more than one square, we can see that $n(\Gamma,z) =1$ for every $z \in \tilde E$ and thus for every $z \in E$, even if $z$ lies on the edge or corner of a square, while $n(\Gamma, z) = 0$ for every $z \in F$. Now any of the remaining segments avoids $\tilde E$ and so must be at least distance $\tfrac{r}{4}$ from $E$. On the other hand, each of these segments adjoins a square which does meet $\tilde E$ and so must be at most distance $\tfrac{r}{2}$ from $E$ and thus at least distance $\tfrac{3r}{2}$ from $F$. 

It is not hard to see by checking cases that no endpoint of a segment in $\Gamma$ which remains after cancellation can be a meeting point of exactly one or exactly three segments. From this, a relatively straightforward argument by induction on the number of segments in $\Gamma$ shows we can express it as a sum of closed curves

\vspace{-.2cm}
\[\Gamma = \Gamma_1 + \Gamma_2 + \cdots \cdots  + \Gamma_n.\] 

Note that the number $n$ of such curves is again uniformly bounded above in terms of $r$ and $R$. 

If $n = 1$ and $\Gamma$ consists of only one curve, then the result follows on using the upper bound on the hyperbolic metric in Lemma 3.1. The main work, then, is in making use of  the curve $\tilde \gamma$ given in the statement to join these potentially disjoint curves together without increasing the hyperbolic length too much. 

Pick $1 \le i \le n$, let $a$ and $b$ be the two points on $\tilde \gamma$ and $\Gamma_i$ respectively which are as close as possible and let $L_i$ denote the line segment joining $a$ and $b$. By replacing those parts of $L_i$ 
(if any) which lie inside one of the other curves $\Gamma_j$, $j \ne i$, with parts of segments of $\Gamma_j$, we can obtain a curve $\tilde L_i$ joining $a$ to $b$ with 
$\dist(\tilde L_i, E) \ge \tfrac{r}{4}$. Note that as the number of squares enclosed by $\Gamma$ is uniformly bounded, so is the number of such replacements we need to carry out and thus the amount of length we add to $L_i$ (in fact, it is not too hard to show we can do this in such a way that the parts of $\Gamma_j$ which are used in this manner can be chosen so as not to overlap).

Next, we examine the distances of points on $\tilde L_i$ to $F$. If $z$ is a point on $\tilde L_i$ within distance $\tfrac{r}{2}$ of $a$, then as $\tilde \gamma$ is at least distance $r$ from $F$, we know that $z$ must be at least distance $\tfrac{r}{2}$ from $F$. Otherwise, $z$ is at least distance $\tfrac{r}{2}$ from $a$ and is either on our original line segment $L_i$ or on $\Gamma$ after we have done the above replacement. Note that in the second case, we have already seen that points on $\Gamma$ are at least distance $\tfrac{3r}{2}$ from $F$, so all we need to check is the first case. 

So suppose now that $z$ is on $L_i$ and at least distance $\tfrac{r}{2}$ from $a$. We claim that $z$ cannot be closer than $\tfrac{r}{2}$ to $F$. Since $a$ is at least $r$ away from $E$ while $b$ is at most $\tfrac{r}{2}$ from $E$, 
$|b - a| \ge \tfrac{r}{2}$. Thus by elementary geometry, if $z$ were within distance less than $\tfrac{r}{2}$ of $F$, we would have a point of $F$ within distance less than $|b - a|$ of $b$. However, since $b$ is within distance $\tfrac{r}{2}$ of $E$, it 
lies in the same complementary component of $\tilde \gamma$ as does $E$. As $\tilde \gamma$ separates $E$ and $F$, it would then follow that we had a point of $\tilde \gamma$ within distance strictly less than $|b-a|$ of $b$ which lies in $\Gamma_i$, which contradicts the fact that $a$ and $b$ are as close as possible. 

Since $\tilde \gamma$ is within distance $R$ of $0$, so must be the set $E$ which it contains as well as the sets $\tilde E$ and hence $\Gamma$. It then follows that each of the line segments $L_i$ above has Euclidean length at most $2R$ and as the length of $\Gamma$ is uniformly bounded, the length of the modified segments $\tilde L_i$ is also uniformly bounded above and this bound is uniform in $r$ and $R$. 

Now using the upper bound on the hyperbolic metric in Lemma 3.1, the hyperbolic length of $\tilde \gamma$, of the curves $\Gamma_i$ and of the curves $\tilde L_i$ above are all bounded uniformly in terms of $r$ and $R$. If we now connect the curves $\Gamma_i$ of $\Gamma$ using the curves $\tilde L_i$ and pieces of $\tilde \gamma$ both of which we traverse in both directions (so that the new curve we obtain is homologous in $U$ to $\Gamma$), then the fact that, as we have already seen, the number $n$ of such curves $\Gamma_i$ is bounded in terms of $r$ and $R$ shows that we can obtain the desired curve $\gamma$ of the statement. $\Box$

{\bf Proof of Theorem 5.1 \hspace{.4cm}}  We will show that $\normV \,\,< \!\infty $ and then appeal to Theorem 4.2. By the bi-equicontinuity of $\calF$, the infinitesimal ratios $\tfrac{|{\mathrm d^\#}f_\alpha'(z)|}{|{\mathrm d^\#}z|}$ (i.e. the expansion of $f_\alpha$ in tangent space) are bounded below at $z = u_\alpha$. By pre- and post-composing with $\tfrac{1}{z}$ which preserves the spherical metric, we can assume that $u_\alpha$ and $v_\alpha = f_\alpha(u_\alpha)$ both lie in the closed unit disc $\overline \D$.  
The fact that $\delta^\#(v_\alpha)$ is bounded below then follows immediately from Bloch's theorem (e.g. \cite{Con} Page 293, Chapter XII, Theorem 1.4).  

For the lower bound on the spherical diameters $\diams (\cbar \setminus V_\alpha)$, since the quantities $\delta_{U_\alpha}^\#(u_\alpha), \delta_{V_\alpha}^\#(v_\alpha)$ are bounded below, we can pre-and post compose with suitable uniformly bi-Lipschitz M\"obius transformations to assume without loss of generality that $U_\alpha, V_\alpha$ are both subdomains of $\C$ and that $u_\alpha = v_\alpha = 0$.
The bi-equicontinuity of $\calF$ then implies that the absolute values of the (Euclidean) derivatives $|f_\alpha'(u_\alpha)|$ are bounded below. The desired conclusion then follows from the corresponding lower bound on the spherical diameters $\diams (\cbar \setminus U_\alpha)$ for the domains of $\calU$ and on applying the Koebe one-quarter 
theorem to a suitable inverse branch of each $f_\alpha$ on the disc ${\mathrm D}(v_\alpha, \delta_{V_\alpha}(v_\alpha))$ where $\delta_{V_\alpha}(v_\alpha)$ is the Euclidean distance from $v_\alpha$ to the boundary $\partial V_\alpha$. 

We now turn to examining the lengths of $\calV$. Let $\gamma_\alpha$ be a meridian in one of the domains $V_\alpha$ of $\calV$ and let $\eta_\alpha$ be a lifting of $\gamma_\alpha$ which lies in $U_\alpha$. Then, since $f_\alpha$ is of degree $\le d$ as a covering map, $\eta_\alpha$ is a closed hyperbolic geodesic in $U_\alpha$ which has length at most $d \ell(\gamma_\alpha)$ as does any meridian which separates the complement of $U_\alpha$ in the same way as $\eta_\alpha$
whence it follows from the fact that $\normU \,< \infty$ that the lengths of $\calV$ must be bounded below. 

We now show that the lengths of $\calV$ are also bounded above. Again, by the uniform lower bound on the quantities $\delta^\#(u_\alpha)$, $\delta^\#(v_\alpha)$, we can assume that the domains $U_\alpha, V_\alpha$ are subdomains of $\C$ with $u_\alpha = v_\alpha = 0$ for each $\alpha$. Again let $\gamma_\alpha$ be a meridian of $V_\alpha$ which separates $\cbar \setminus V_\alpha$ into disjoint non-empty closed sets $E_\alpha$, $F_\alpha$ with $\infty \in F_\alpha$. If we once again let $\eta_\alpha$ be a lifting of $\gamma_\alpha$, then $\eta_\alpha$ is a simple closed geodesic in $U_\alpha$ and by Theorem 1.6 we can find a meridian $\tilde \eta_\alpha$ in the homology class of $\eta_\alpha$. Using Cauchy's theorem, if we let $\tilde \gamma_\alpha$ be the curve $f_\alpha(\tilde \eta_\alpha)$ (which may possibly be traversed more than once), then one can check that the winding number of $\tilde \eta_\alpha$ about points of $E_\alpha$ is non-zero while the winding number about points of $F_\alpha$ is zero. 
Hence $\tilde \gamma_\alpha$ separates $E_\alpha$ and $F_\alpha$.  

Since $\normU \, < \infty$, it follows from the remarks after Definition 4.5 that the hyperbolic length of $\tilde \eta_\alpha$ 
and the distance of this meridian to the basepoint $u_\alpha$ are uniformly bounded above. Hence by the Schwarz lemma for the hyperbolic metric, the hyperbolic lengths of the curves $\tilde \gamma_\alpha$ and their distances from the base points $v_\alpha$ in $V_\alpha$ are also uniformly bounded above. As $v_\alpha = 0$ and $\delta^\#(v_\alpha)$ is bounded below, it follows from Corollary 4.1 that we can find $R > 0$ such that the curves $\tilde \gamma_\alpha$ all lie in $\overline{{\mathrm D}(0,R)}$ and the hyperbolic lengths of these curves are all $\le R$. It also follows from Corollary 4.1 that we can find $r >0$ such that for every $\alpha$, $\dist (\tilde \gamma_\alpha, \partial V_\alpha) \ge r$. 

We thus have satisfied the hypotheses of the above lemma for the domain $V_\alpha$, the sets $E_\alpha$, $F_\alpha$ and the curve $\tilde \gamma_\alpha$. Applying the lemma then gives us a curve in the homology class of $\gamma_\alpha$ whose hyperbolic length in $V_\alpha$ is uniformly bounded above in $\alpha$. The desired uniform upper bound on the length of $\gamma_\alpha$ is then immediate in view of Theorem 1.6.  

Finally, we examine the distances of $\calV$. As observed above, the hyperbolic distance from $u_\alpha$ to $\tilde \eta_\alpha$ is uniformly bounded above and by \cite{Com1} Lemma 2.2, $\eta_\alpha$ and $\tilde \eta_\alpha$ must intersect (which includes the possibility that they are the same curve). Again using $\normU \,< \infty$, the hyperbolic distance from $u_\alpha$ to $\eta_\alpha$ is uniformly bounded above and by the Schwarz lemma the hyperbolic distance from $v_\alpha$ to $\gamma_\alpha$ is then also uniformly bounded above. $\Box$

\section{Bounded Containment}

Instead of being bounded (in $\cbar$), a family of pointed domains can be bounded within another bounded family. 

\begin{definition}
Let $\Upup$ and $\Uu$ be pointed domains of connectivity $n'$ and $n$ respectively with $U' \subset U$. We say that $U'$ is {\rm bounded above and below} or just {\rm bounded} in $U$ with constant $K \ge 1$ if: 
\vspace{-.1cm}

\begin{enumerate}
\item $U'$ is a subset of $U$ which lies within hyperbolic distance at most $K$ about $u$ in $U$;

\vspace{.2cm}
\item \[ \delta^\#_{U'} (u') \ge \frac{1}{K}  \delta^\#_{U} (u');\]

\item \[\calL(\Upup) \le K;\]

\item \[\calD(\Upup) \le K.\]

\end{enumerate}

If $\,\calU = \{\Uaua \}_{\alpha \in A}$ is a family of $n$-connected pointed domains, we say that another family $\calUp = \{\Uauap \}_{\alpha \in A}$ of $n'$-connected pointed domains is {\rm bounded above and below} or just {\rm bounded in $\calU$} if we can find $K \ge 1$ such that for every $\alpha$, $U'_\alpha$ is bounded in $U_\alpha$ with this constant $K$. 
In this case we write $\pt \sqsubset \calUp
\sqsubset \calU$.
 \end{definition}

Note the similarities with the \Car bound as given in Definition 4.6. Once again, if the domains of $\calUp$ are simply connected, then conditions \emph{3.} and \emph{4.} are vacuously true. 

We are also interested in boundedness for the degenerate case of a family of simple closed curves, such as meridians of a family of domains (this is useful in considering quasiconformal gluing problems, for example when gluing the boundaries of certain annuli in order to prove a non-autonomous version of the Sullivan straightening theorem e.g. \cite{Com4}).

\begin{definition} Let $\Gamma = \{(\gamma_\alpha, z_\alpha)\}_{\alpha \in A}$ be a family of pointed simple smooth (i.e. $C^1$) closed curves where for each $\alpha$, $z_\alpha$ is a point on $\gamma_\alpha$ and $\phi_\alpha : \T \mapsto \cbar$ is a parametrization of $\gamma_\alpha$. 

Let $\calU$ be a family of pointed domains of the same finite connectivity. We say $\Gamma$ is {\rm bounded above and below} 
or just {\rm bounded} in $\calU$ with constant $K \ge 1$ if 
\begin{enumerate}
\item For each $\alpha$, $\gamma_\alpha$ is a subset of $U_\alpha$ which lies within hyperbolic distance at most $K$ about $u_\alpha$ in $U_\alpha$; 

\vspace{.2cm}
\item The mappings $\phi_\alpha$ can be chosen so that the resulting family 
$\Phi$ is bi-equicontinuous on $\T$ in the sense that we can find $K \ge 1$ such that 

\[ \frac{1}{K} \le \frac{|\phi_\alpha^\#(z)|}{\delta_{U_\alpha}^\#(z_\alpha)} \le K, \qquad z \in \T,\: \alpha \in A.\]

\end{enumerate}

In this case, we write $\pt \sqsubset \Gamma \sqsubset \calU$.

We say $\Gamma$ is {\rm bounded above and below} or just {\rm bounded} in $\cbar$ with constant $K \ge 1$ if the mappings $\phi_\alpha$ can be chosen so that the resulting family 
$\Phi$ is bi-equicontinuous on $\T$ in the sense that we can find $K \ge 1$ such that 

\vspace{-.2cm}
\[ \frac{1}{K} \le |\phi_\alpha^\#(z)| \le K, \qquad z \in \T, \: \alpha \in A.\]

In this case, we write $\pt \sqsubset \Gamma \sqsubset \cbar$.

\end{definition}

Using Theorem 4.2 and Proposition 5.2, we immediately have the following.

\begin{proposition} Let $\calU = \{\Uaua\}_{\alpha \in A}$ be a bounded family of $n$-connected domains, for each $\alpha$ let $\gamma_\alpha$ be a meridian of $U_\alpha$ and let $\Gamma$ be the resulting family. Then $\pt \sqsubset \Gamma \sqsubset \calU$.
\end{proposition}

From now on unless otherwise specified we will use 
$\calUp, \calU, \calVp, \calV$ and $\calW$ to refer to families of multiply connected domains where all the domains in a given family have the same connectivity.

The next two results are natural consequences of the notion of boundedness.

\begin{theorem} [Transitivity of Boundedness]
If $\pt \sqsubset \calU \sqsubset \calV \sqsubset \calW$ where $\calW$ is either bounded or all the domains of $\calW$ are $\cbar$, then 
$\pt \sqsubset \calU \sqsubset \calW$. This includes the degenerate case where $\calU$ is a family of simple closed curves.
\end{theorem}

\proof Consider first the case when $\calU$ is a family of pointed domains and $\calW$ is bounded. The only thing which needs to be checked is {\emph 2.} and by \emph{1.} and \emph{2.} for the boundedness of $\calU$ in $\calV$ and of $\calV$ in $\calW$, Corollary 4.1 and the Schwarz lemma for the hyperbolic metric
\vspace{.25cm}
\[\delta^\#_{U_\alpha}(u_\alpha) \sim \delta^\#_{V_\alpha}(u_\alpha) \sim \delta^\#_{V_\alpha}(v_\alpha) \sim \delta^\#_{W_\alpha}(v_\alpha) \sim \delta^\#_{W_\alpha}(w_\alpha) \sim \delta^\#_{W_\alpha}(u_\alpha)\]

where as usual the subscripts denote the domains with respect to which the distance to the boundary is taken. 

In the case where $\calU$ is a family of hyperbolic domains and all the domains of $\calW$ are $\cbar$, we shall show that 
$\normU \,< \infty$ and then appeal to Theorem 4.2.  By {\emph 1.} and {\emph 2.} for the boundeness of $\calU$ in $\calV$ combined with Corollary 4.1 and the bound for $\delta^\#_{V_\alpha}(v_\alpha)$ arising from the fact that the boundedness of $\calV$ implies that $\normV \,< \infty$
\vspace{.25cm}
\[\delta^\#_{U_\alpha}(u_\alpha) \sim \delta^\#_{V_\alpha}(u_\alpha) \sim \delta^\#_{V_\alpha}(v_\alpha) \sim 1\]

whence we have the correct bounds on the spherical distances from each $u_\alpha$ to the boundary of $U_\alpha$. As $U_\alpha \subset V_\alpha$ for each $\alpha$, the spherical diameters $\diams (\cbar \setminus U_\alpha)$ are bounded below and so we have taken care of the first term in the \Car bound.
Finally, the bounds on the lengths of $U_\alpha$ and on the hyperbolic distances in $U_\alpha$ from the basepoints $u_\alpha$ to the meridians of a maximally close system for $U_\alpha$ follow directly 
from \emph{ 3.} and \emph{4.} for the boundeness of $\calU$ in $\calV$.

If $\calU$ is a family of curves and $\calW$ is bounded, the result follows easily from Corollary 4.1 and the fact that if $z$ is a point on the curve $U_\alpha$, then, similarly to above
\vspace{.25cm}
\[\delta^\#_{V_\alpha}(z) \sim \delta^\#_{V_\alpha}(v_\alpha) \sim \delta^\#_{W_\alpha}(v_\alpha) \sim \delta^\#_{W_\alpha}(w_\alpha) \sim \delta^\#_{W_\alpha}(u_\alpha) \sim \delta^\#_{W_\alpha}(z).\]

For the final possibility, if $\calU$ is a family of curves and all the domains of $\calW$ are $\cbar$, if  $z \in U_\alpha$, by the boundedness of $\calU$ in $\calV$, the boundedness of $\calV$ and again Corollary 4.1
\vspace{.25cm}
\[\delta^\#_{V_\alpha}(z) \sim \delta^\#_{V_\alpha}(u_\alpha) \sim \delta^\#_{V_\alpha}(v_\alpha) \sim 1.\]

This shows the curves $U_\alpha$ also give a bi-equicontinuous family with respect to the spherical metric and with this the proof is now complete. $\Box$

\begin{proposition} [Bi-Equicontinuity and Bounded Containment]
Let\\ $\calU$, $\calUp$ be two families with $\pt \sqsubset \calUp \sqsubset \calU \sqsubset \cbar$ and let $\calF$ be a family of 
covering maps (for both the domains of $\calU$ and $\calUp$) of uniformly bounded degree with 
$\calF \bowtie \,\calU$. Then if $\calV, \calVp$ are the corresponding image families, we have $\pt \sqsubset \calVp \sqsubset \calV \sqsubset \cbar$ and the constant for this boundedness is uniform with respect to the boundedness of $\calUp$ in $\calU$ and the degree bound for the mappings of $\calF$. This includes the degenerate case where $\calF$ is univalent and $\calUp$ and $\calVp$ are families of simple closed curves \end{proposition}  

\proof For the non-degenerate case where we have families of pointed domains, the bound on the hyperbolic distance of points in $V'_\alpha$ from the basepoint $v_\alpha$ for the larger domain $V_\alpha$ follows from the Schwarz lemma. 

For the estimate on $\delta^\#_{V'_\alpha}(v'_\alpha)$ we examine the inverse branches of the mappings $f_\alpha$. Note that as before we can postcompose with $\tfrac{1}{z}$ which preserves spherical distances so that we can assume that the points $v'_\alpha$ lie in the closed unit disc and it will thus suffice to look at the Euclidean distances $\delta_{V'_\alpha}(v'_\alpha)$, $\delta_{V_\alpha}(v'_\alpha)$ to the boundaries of $V'_\alpha$, $V_\alpha$ respectively.

Now apply the Koebe one-quarter theorem on the discs ${\mathrm D}(v'_\alpha, \delta_{V_\alpha}(v'_\alpha))$ and then the distortion theorem for univalent mappings (e.g. \cite{CG} Page 3, Theorem 1.6) to inverse branches of $f_\alpha$ on the discs ${\mathrm D}(v'_\alpha, \delta_{V'_\alpha}(v'_\alpha))$. It then follows that $\delta_{V'_\alpha}(v'_\alpha)$ cannot be too small compared to $\delta_{V_\alpha}(v'_\alpha)$ since otherwise this would force $\delta_{U'_\alpha}(u'_\alpha)$ to be very small compared to $\delta_{U_\alpha}(u'_\alpha)$ which would contradict the boundedness of $\calUp$ in $\calU$. 

Since $\calUp \sqsubset \calU$ and $\calU$ is bounded, by Theorem 6.1, $\calUp$ is also bounded. The required bounds on the distances and lengths of both $\calV$ and $\calVp$ then follow from Theorem 5.1. Finally,  the proof in the case of families of simple closed curves is straightforward in view of Theorem 4.1. $\Box$

The following is a more or less direct consequence of the definition and Corollary 4.1. 

\begin{corollary}
Let $\calV$ be a bounded family of $n$-connected domains for some $n \ge 1$ and let $\calU$ be either a family of $m$-connected domains for some $m \ge 1$ or a family of curves. Then if $\pt \sqsubset \calU \sqsubset \calV \sqsubset \cbar$, there exists $\delta >0$ depending only on $m$, $n$ and $K$ or $n$ and $K$ as appropriate such 
that, for every $\alpha \in A$, $V_\alpha$ contains a (spherical) 
$\delta$-neighbourhood of $U_\alpha$. 
\end{corollary}

Sometimes it is useful to know whether we can get a new bounded family if we remove the closures of the domains for one bounded family from another and then make a suitable choice of base points. If the two original families are families of discs and one is bounded in the other, the answer is yes.

\begin{theorem}
Let $\calU = \{\Uaua\}_{\alpha \in A}$ and $\calV = \{\Vava\}$ be two families of pointed discs with  
$\pt \sqsubset \calU \sqsubset \calV \sqsubset \cbar$. For each $\alpha$, let $A_\alpha$ be the conformal annulus 
$V_\alpha \backslash \overline {U_\alpha}$, let $\gamma_\alpha$ be the equator of $A_\alpha$ and let $a_\alpha$ be a point on $\gamma_\alpha$. Then the family $\calA = \{\Aaaa\}_{\alpha \in A}$ is a bounded family of pointed annuli.
\end{theorem}

\proof For each $\alpha$, let $\pi_\alpha$ be the unique normalized (inverse) Riemann
mapping which sends $\D$ to $V_\alpha$. By Proposition 5.2 these mappings form a family $\Pi$
which is univalent and bi-equicontinuous on compact subsets of 
$\D$. Now for each $\alpha$ set
$\tilde U_\alpha = \pi_\alpha^{\circ -1}(U_\alpha)$, $\tilde A_\alpha = \D \backslash \overline {\tilde U_\alpha}$. 
It follows from the fact that $\pt \sqsubset \calU \sqsubset \calV \sqsubset \cbar$ that the moduli and hence the lengths of the annuli $\tilde A_\alpha$ are uniformly bounded above and below. Thus, by conformal invariance, the same is true of the moduli and lengths of the annuli $A_\alpha$. If for each $\alpha$, we now pick an arbitrary point $a_\alpha$ on the equator $\gamma_\alpha$ of 
$A_\alpha$, then it follows that the second and third terms in the \Car bound of each $A_\alpha$ will be uniformly bounded. 

To finish, we need to establish a uniform bound on the first term, namely on 
$\delta^\#_{A_\alpha} (a_\alpha)$ and $\diams (\cbar \setminus A_\alpha)$. It follows from the above bounds on the moduli and from the usual estimates on the hyperbolic metric in Lemma 3.1 that the equators $\tilde \gamma_\alpha$ of the annuli 
$\tilde A_\alpha$ must be bounded away from the boundary curves in terms of Euclidean distance. Since one of these boundary curves is simply the unit circle, it follows that they also must all lie inside
a disc of some bounded hyperbolic radius about $0$ in $\D$. Again we can apply $\tfrac{1}{z}$ if needed so as to assume that the points $a_\alpha$ lie in $\overline \D$. This allows us to apply 
the Koebe one-quarter theorem and the bi-equicontinuity of $\Pi$ to deduce that the equators of the annuli $A_\alpha$ will also be 
uniformly bounded away from their boundary curves in terms of Euclidean and hence spherical distances. Since the unbounded complementary component of each annulus $A_\alpha$ is simply $\cbar \setminus V_\alpha$, the uniform lower bound on the spherical diameters of the complements $\cbar \setminus A_\alpha$ follows from Theorem 4.2 and the corresponding lower bound on the spherical diameters of the sets $\cbar \setminus V_\alpha$. $\Box$ 

The reader might wonder if it is necessary to restrict ourselves to families of discs in the above statement. However, when one tries to consider domains of higher connectivity, problems can arise. Consider, for example the case in Figure 6 where for each $m \ge 1$
\[A_m = {\mathrm A}(0,2,5), \qquad D_m = \{z: 3  < |z| < 4, 1/m < {\rm Arg}\, z < 2\pi - 1/m\}.\] 

It is then easy to see that the pointed annuli $(A_m, -3.5)$ and the pointed discs $(D_m, -3.5)$ both give rise to bounded families which we will call 
$\calA$ and $\calD$ respectively. It is then also obvious to see that 
$\pt \sqsubset \calD \sqsubset \calA \sqsubset \cbar$.

\scalebox{1.05}{\includegraphics{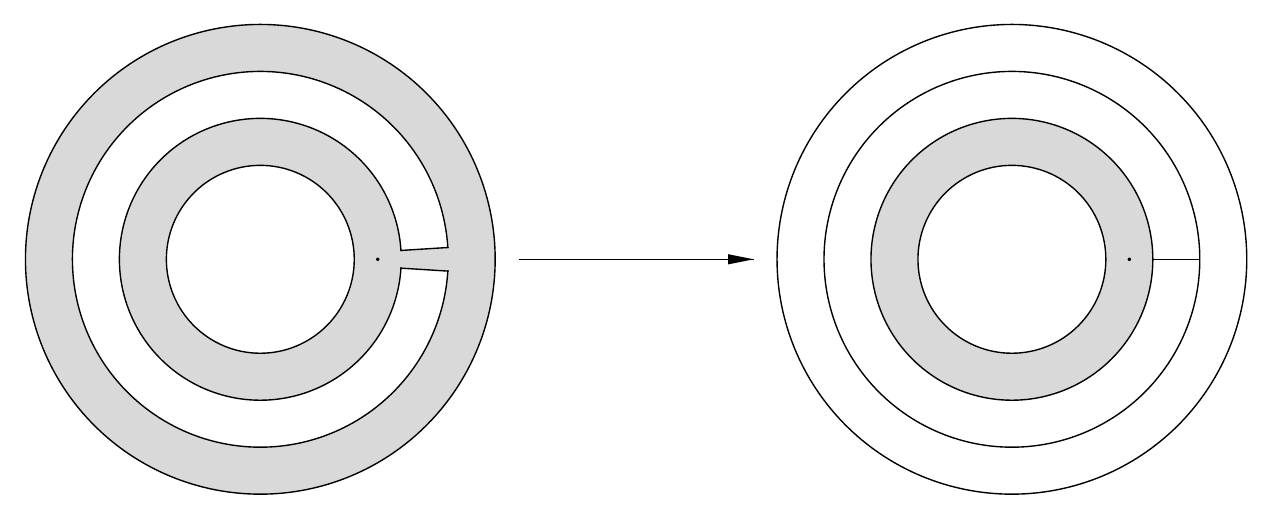}}
\unitlength1cm
\begin{picture}(0.01,0.01)
  \put(-11.57,-.2){\footnotesize$(U_m,u_m)$  }
  \put(-7.47,3){\footnotesize$m \to \infty$}
  \put(-3.27,-.2){\footnotesize$(U,u)$  }
\end{picture}

\vspace{.7cm}
Now remove the closure of the pointed disc $D_m$ from the annulus $A_m$ and let $U_m$ be the resulting triply connected region. Now let $\{u_m\}_{m=1}^\infty$ be a convergent sequence where for each $m$, $u_m \in U_m$ and let $\calU$ be the resulting family of pointed domains. Depending on our choice of base points $u_m$, the \Car limit for the pointed domains $\Umum$ will be either a point or doubly connected, essentially because the discs $D_m$ will `close off' at least one of the complementary components of $U_m$ as seen from the base point $u_m$. We thus have a family of triply connected domains where the only possibility for a limit is either a point or doubly connected and so $\calU$ cannot be bounded.

To conclude this paper, as an application we solve an extremal problem. Suppose once more we have a multiply connected domain $U$ whose complement can be expressed as the union of two disjoint non-empty closed subsets $E$ and $F$ neither of which is a point. We know from Theorem 1.6 that we can find a meridian which separates $E$ and $F$ and hence by the collar lemma there is a conformal annulus of some definite but bounded modulus which also separates $E$ and $F$. 

Since neither $E$ nor $F$ is a point, any annulus which separates them must have uniformly bounded modulus. Now let $M < \infty$ be the supremum over the moduli of all such annuli and take a sequence of pointed annuli $(A_m, a_m)$ with base points on their equators and whose moduli tend to $M$ as $m$ tends to infinity. The diameters of the complementary components of these annuli must be bounded below since they contain $E$ and $F$ and using Corollary 4.1 (or Lemma 3.1), it follows easily that the spherical distances $\delta^\#_{A_m}(z)$ of points on their equators (including the base points) to the boundary will be uniformly bounded below and above. By \emph{4.} $\Longrightarrow$ \emph{1.} of Theorem 4.2, this gives us a bounded family and if we now take a subsequence if necessary so that the base points $a_m$ converge to some limit $a$ and the complementary components converge in the Hausdorff topology, then the domains $(A_m, a_m)$ converge in the \Car topology to a pointed annulus $\Aa$. By \emph{2.} of Theorem 1.10, the modulus of $A$ must then be $M$. 

We still need to show that $A$ separates $E$ and $F$. Note that by ii) of \Car convergence, $A \cap (E \cup F) = \emptyset$. If we let $\gamma_m$ be the equator of each $A_m$ and $\gamma$ the equator of $A$, then by Theorem 1.10 the curves $\gamma_m$ converge uniformly to $\gamma$ and are thus eventually homotopic in $U$ to $\gamma$. 
Since each curve $\gamma_m$ separates $E$ and $F$, it then follows easily using winding numbers that $\gamma$ and thus $A$ separates $E$ and $F$. We have thus proved the following.

\begin{theorem} Let $U$ be multiply connected and suppose $\cbar \setminus U = E  \cup F$ where $E$ and $F$ are closed and disjoint and neither set is a point. Then we can find a conformal annulus $A$ in $U$ which separates $E$ and $F$ for which 
\[ \bmod A = \sup_B\bmod B\]
where the supremum is taken over all conformal annuli $B$ in $U$ which separate $E$ and $F$.
\end{theorem}

Essentially this is a variant of the extremal problems considered by Gr\"otzsch, Teichm\"uller and Mori (e.g. \cite{Ahl2, McM} and more recently in the paper of Herron, Liu and Minda \cite{HLM}. We remark that it is also possible to give a somewhat longer but more elementary proof of this result using Theorem 1.2 and we are indebted to Adam Epstein for pointing this out. 

Although the maximum possible modulus is attained, there is still the question as to whether such an annulus is unique and the following counterexample shows this is not in general true. Let $E = \overline {\mathrm D}(-3,1) \cup \overline {\mathrm D}(3,1)$, let $F = \overline \D \cup (\cbar \setminus {\mathrm D}(0,5))$ and let $U$ be the quadruply connected domain $\cbar \setminus (E \cup F)$. 

From above there is a conformal annulus $A$ of maximum possible modulus in $U$ which separates $E$ and $F$ and let us suppose for the sake of contradiction that $A$ is unique. Now $U$ is symmetric under the transformation $z \mapsto -z$ and, by the assumed uniqueness, so is $A$. If we then let $\gamma$ be the equator of $A$, then by the symmetry of $A$ and the uniqueness of the equator of a conformal annulus, $\gamma$ must also in turn be symmetric under this transformation. Using winding numbers, it follows that $\gamma$ encloses either just $\overline \D$ or 
$\overline \D \cup \overline {\mathrm D}(-3,1) \cup \overline {\mathrm D}(3,1)$. In neither case, does $\gamma$ separate $E$ and $F$ which is clearly impossible as $A$ must separate these sets and, with this contradiction, we have what we want.


\end{document}